\documentclass{birkjour}
\usepackage{amsmath,amssymb,isomath}
\usepackage{txfonts,fourier}
\usepackage{delarray}
\usepackage{xcolor,graphicx}
\usepackage{pgfplots}
\usepackage{tikz,pgf}
\usepackage[colorlinks=true,urlcolor=blue,
citecolor=red,linkcolor=blue,linktocpage,pdfpagelabels,
bookmarksnumbered,bookmarksopen]{hyperref}

\newcommand{\eps}{\varepsilon}

\newcommand{\om}{\omega}

\renewcommand{\th}{\theta}

\newcommand{\cU}{\mathcal{U}}
\newcommand{\cS}{\mathcal{S}}
\newcommand{\cZ}{\mathcal{Z}}
  \newcommand{\z}{{\zeta}}

\newcommand\enne{\varmathbb{N}}

\newcommand\erre{\varmathbb{R}}

\newcommand{\pa}{\partial}

\newcommand{\ov}{\overline}
\newcommand{\dyle}{\displaystyle}
\newcommand{\dimo}{{\bf Proof.}\quad}
\newcommand{\finedim}{{\unskip\nobreak\hfil\penalty50
  \hskip2em\hbox{}\nobreak\hfil\mbox{\rule{1ex}{1ex} \qquad}
  \parfillskip=0pt\finalhyphendemerits=0\par\medskip}}
  \newcommand{\re}{{\rm Re\,}}
    \newcommand{\im}{{\rm Im\,}}

\newtheorem{theorem}{Theorem}[section]
\newtheorem{remark}[theorem]{Remark}
\newtheorem{lemma}[theorem]{Lemma}
\newtheorem{proposition}[theorem]{Proposition}
\newtheorem{corollary}[theorem]{Corollary}
\numberwithin{equation}{section}
\newcommand{\bteo}{\begin{theorem}}
\newcommand{\eteo}{\end{theorem}}
\newcommand{\bprop}{\begin{proposition}}
\newcommand{\eprop}{\end{proposition}}
\newcommand{\blem}{\begin{lemma}}
\newcommand{\elem}{\end{lemma}}
\newcommand{\boss}{\begin{remark}\rm}
\newcommand{\eoss}{\end{remark}}
\newcommand{\beq}{\begin{equation}}
\newcommand{\eeq}{\end{equation}}

\begin{document}

\title[On four species strongly competing systems]{On the limit 
configuration of four species \\
strongly competing systems}

\author{Flavia Lanzara}
\address{\rm Dipartimento di Matematica {\it G. Castelnuovo}\\
{\it Sapienza} Universit\`a di Roma,
piazzale Aldo Moro 5, 00185 Roma}
\email{lanzara@mat.uniroma1.it}

\author{Eugenio Montefusco}
\address{\rm Dipartimento di Matematica {\it G. Castelnuovo}\\
{\it Sapienza} Universit\`a di Roma,
piazzale Aldo Moro 5, 00185 Roma}
\email{montefus@mat.uniroma1.it}



\begin{abstract}
We analysed some qualitative properties of the limit configuration 
of the solutions of a reaction-diffusion system of four competing 
species as the competition rate tends to infinity.  
Large interaction induces the spatial segregation of the species and only 
two limit configurations are possible: either there is a point where four 
species concur, a 4-point, or there are two points where only three 
species concur. We characterized, for a given datum, the possible 4-point 
configuration by means of the solution of a Dirichlet problem for the 
Laplace equation. 
\end{abstract}
\maketitle

\section{Introduction and setting of the problem}

In ecology, the competitive exclusion principle (also known as 
Gause's law) is the proposition that two (or more) species 
competing exactly for the same limiting resource cannot coexist 
at constant population values: if one specie has even the slightest 
advantage over another, the one with the advantage will dominate 
in the long term. This leads either to the extinction of the weaker 
competitors or to an evolutionary or behavioral shift toward a different 
ecological niche. The principle has been synthetized in the statement 
{\it complete competitors cannot coexist}.
According to the competitive exclusion principle, many 
competing species cannot coexist under very strong competition,
but when spatial movements are permitted more than one species can 
coexist thanks to the segregation of their habitats. For a theoretical
discussion and some experimental results see \cite{gau,ha}.

From a mathematical viewpoint the determination of the 
configuration of the habitat segregation for some populations is an 
interesting problem which can be modelled by an optimal (in a suitable 
sense) partition of a domain, for example in the papers 
\cite{ctv1,ctv2,ctv3,ctv4,tt12} the problem is studied modelling the 
interspecies competition with a large interaction term in an elliptic 
system of partial differential equations inspired by classical models
in populations dynamics. In \cite{bz,dwz} the problem is modelled as a 
Cauchy problem for a parabolic system of semilinear partial differential 
equations describing the dynamics of the densities of different species. 
Starting from the quoted papers we want to tackle the segregation problem 
of 4 species in a planar region. Note that in the evolution works, in
particular \cite{dwz}, it is proved that some populations can vanish under
the competition of other species, moreover in \cite{bz} the authors are 
able to estimate the number of the long-term surviving populations.

Let \(D\subseteq\erre^n\) be an open bounded, simply connected  
domain with smooth boundary $\pa D$ and consider the 
competition-diffusion system of \(k\) differential equations
\begin{equation}\label{diffeqk}
\begin{array}\{{cl}.
-\Delta u_i(x)=-\mu u_i (x)\dyle\sum_{j\neq i} u_j (x)& \hbox{ in }D\,, \\
u_i(x)\geq 0 & \hbox{ in }D\,, \\
u_i(x)=\phi_i(x) & \hbox{ on }\pa D\,.
\end{array} \qquad i=1,...,k
\end{equation}
This system governs the steady states of \(k\) competing species 
coexisting in the same domain \(D\). Any \(u_i\) represents a population 
density and the parameter \(\mu\) determines the interaction 
strength between the populations. \\
The datum  \(\Phi=(\phi_1,...,\phi_k)\) is called {\bf admissible} if
any function \(\phi_i\in H^{1/2}(\partial D)\) (for \(i=1,...,k\)) is 
nonnegative, positive in a nonempty arc, and such that 
\(\phi_i\cdot\phi_j=0\) a.e. in \(\partial D\) for \(i\neq j\). \\ 
Moreover, if \(n=2\), we will assume that \(\phi_i\in W^{1,\infty}(\pa D)\), 
that the sets \(\{\phi_i>0\}\) are nonempty, open connected arcs and 
that the function \(\sum_{i=1}^k\phi_i\) vanishes exactly in \(k\) points 
of \(\partial D\) (the endpoints of the \(\phi_i\)'s supports).\\
If \(\Phi\) is admissible, the existence of positive solutions of 
\eqref{diffeqk} for any positive \(\mu\) is proved in \cite{ctv2} using 
Leray-Schauder degree theory. The uniqueness is proved in \cite{wz}, 
using the sub- and super-solution method. \\
Let define the class of segregated densities
\[
\cU = \begin{array}\{{rl}\}
U=(u_1,...,u_k)\in (H^1(D))^k \hbox{ : }  & u_i=\phi_i \hbox{ on }\pa D \\
& u_i\geq 0 \hbox{ in } D \\
 & u_i\cdot u_j=0 \hbox{ for \(i\neq j\) a.e. in } D \\
\end{array}
\]
and the class
\[
\cS=\begin{array}\{{rl}\}
U=(u_1,...,u_k)\in \cU\hbox{ : } 
 & -\Delta u_i\leq 0 \hbox{ in } D \\
 & -\Delta \big(u_i-\sum_{j\neq i} u_j \big)\geq 0 \hbox{ in } D \\
\end{array}\,.
\]
Let \(U^{(\mu)}=(u_{1,\mu},...,u_{k,\mu})\) be the solution of 
\eqref{diffeqk} for every \(\mu>0\). In \cite{ctv2} it is proved 
that there exists \(\overline{U}=\{\bar u_1,...,\bar u_k\}\in\cU\) such 
that, up to subsequences, \(u_{i,\mu}\longrightarrow\bar{u}_i\) in
\(H^1(D)\) and \(\cS\) contains all the asymptotic limits of \eqref{diffeqk} 
that is \(\overline{U}\in \cS\). 

The uniqueness of the limit solution of \eqref{diffeqk} as 
\(\mu\to+\infty\) was proved in \cite{ctv2} in the case \(k=2\) 
and in \cite{ctv4} in the case of \(k=3\) and in dimension \(2\). 
Specifically, the authors prove that the class \(\cS\) consists 
of one element.
In \cite{wz}  it is proved that  \(\cS\) consists of one element 
also in the case of arbitrary dimension and arbitrary number of species. 
A different proof of uniqueness of the limit configuration, based 
on the maximum principle and on the qualitative properties of the 
elements of \(\cS\), is given in \cite{ab}.  

Now consider the energy functional associated with \(k\) species
\[
E(U)=\sum_{i=1}^k \int_D |\nabla u_i(x)|^2dx
\]
in the class of all possible segregated states  subject to some 
boundary and positivity
conditions.  It is shown in \cite{ctv3} that the problem
\begin{equation}\label{minimo}
\hbox{find}\quad U\in\cU
\quad\hbox{such that}\quad E(U)=\min_{V\in\mathcal{U}} E(V)
\end{equation}
admits a unique solution, which belongs to $\cS\subseteq\cU$. 
It follows that the unique minimal energy configuration shares 
with the limit states of \eqref{diffeqk} the common property of 
belonging to \(\cS\).
Specifically,  the element of \(\cS\) is the critical point 
(in the weak sense) of the energy functional \(E(U)\) among 
all segregated states, subject to  the same boundary and positivity
conditions.

The aim of this note is to study the case of four species in a 
planar smooth domain (i.e. the case \(k=4\) and \(n=2\) in \eqref{diffeqk}), 
in particular we are interested in the description and some qualitative 
properties of the limiting configuration. 

\section{Some known results}

Suppose that \(D\) is a simply connected domain in \(\erre^2\). 
Due to the conformal invariance of the problem, with no loss 
of generality we can assume
\[
D=B(O,1)=\{x=(x_1,x_2) \in\erre^2: |x|<1\}
\]
and consider the class
\[
\cS=\begin{array}\{{rl}\}
&U=(u_1,u_2,u_3,u_4)\in (H^1(D))^4 \hbox{ : } 
u_i\geq 0 \hbox{ in } D ,\>  u_i=\phi_i \hbox{ on }\pa D\\
&u_i\cdot u_j=0 \hbox{ for } i\neq j,  -\Delta u_i\leq 0,   
-\Delta \big(u_i-\sum_{j\neq i} u_j \big)\geq 0 \hbox{ in } D 
\end{array}.
\]
The study of \(\cS\) provides the understanding of the segregated 
states of \(4\) species induced by strong competition.
If  $\Phi=(\phi_1,\phi_2,\phi_3,\phi_4)$ is an admissible 
datum then \(\phi_i, i=1,...,4,\) are positive in their supports, 
the sets  $\{\phi_i>0\}\subset \partial D$ are open 
connected arcs  and $\phi=\sum_{i=1}^4 \phi_i$ vanishes at exactly 
four points of $\pa D$, the endpoints  $p_1,p_2,p_3,p_4$ in anticlockwise order. 

In the following we will denote by \(U\) both the quadruple 
\((u_1,u_2,u_3,u_4)\) and the function \(\sum_{i=1}^4 u_i\). 
For any $U\in\cS$ define the {\it nodal regions}
\[
\om_i=\{x\in D:u_i(x)>0\} \qquad i=1,...,4\,,
\]
the {\it multiplicity of a point} $x \in\ov{D}$ with respect to $U$:
\[
m(x)=\#\{i: |\om_i\cap B_r(x)|>0 \quad \forall r>0\}
\]
where $B_r(x)=\{p\in\erre^2: |p-x|<r\}$ and the {\it interfaces} 
between two densities
\[
\Gamma_{ij}=\pa \om_i\cap \pa \om_j \cap \{x\in D: m(x)=2\}\,.
\]
We say that $\om_i$ and $\om_j$ are adjacent if $\Gamma_{ij}\neq\emptyset$.
Let us  summarize the basic properties of the elements $U\in\cS$:
\begin{itemize}
\item[(s1)]
$\sum_{i=1}^4 u_i\in W^{1,\infty}(\ov{D})$ (\cite[Theorem 8.4]
{ctv3}). It follows that $u_i\in C(\ov{D})$, $\om_i$ is open 
and \(x\in\om_i\) implies $m(x)=1$;
\item[(s2)] each $\om_i$ is connected and 
each $\Gamma_{ij}$ is either empty or a connected arc starting from the 
boundary  (\cite[Remark 2.1]{ctv4});
\item[(s3)] $u_i$ is harmonic in $\om_i$; $u_i-u_j$ is harmonic in 
$D\setminus \cup_{h\neq i,j}\overline \om_h$ \cite[Proposition 6.3]{ctv3};
\item[(s4)] if $x_0\in D$ such that \(m(x_0)=2\) then (\cite[Remark 6.4 ]{ctv3})
\[
\lim_{\om_i\ni y\to x_0} \nabla u_i(y)
=-\lim_{\om_j\ni y\to x_0} \nabla u_j(y) ;
\]
\item[(s5)] if $x_0\in D$ such that \(m(x_0)=2\) then \(\nabla U(x_0)\neq 0\) 
and the set  $\{x:m(x)=2\}$ is locally a $C^1$-curve 
through $x_0$ ending either at points with higher multiplicity, 
or at the boundary $\pa D$  \cite[Lemma 9.4]{ctv3};
\item[(s6)] if $x_0\in D$ such that $m(x_0)\geq 3$ then $|\nabla U(x)|\to 0$, 
as $x\to x_0$ \cite[Theorem 9.3]{ctv3};
\item[(s7)] the set $\{x: m(x)\geq 3\}$ consists of a finite 
number of points \cite[Lemma 9.11]{ctv3};
\item[(s8)] if $x_0\in D$ with $m(x_0)=h\geq 3$ then there exists 
\(\th\in (-\pi,\pi]\) such that
\begin{equation}\label{local}
U(r,\th)=r^{h/2} \left| \cos\left( \frac{h}{2} (\th +\th_0)\right)
\right|+o(r^{h/2})
\end{equation}
as $r\to 0$, where $(r,\th)$ is a system of polar coordinates around 
\(x_0\) \cite[Theorem  9.6]{ctv3}.
\end{itemize}

\boss The asymptotic formula in (s8) describes the behavior of \(U\) 
in a neighborhood of a multiple point in \(D\). As a consequence, at 
multiple point \(U\in\cS\) shares the angle in equal parts. 
This property does not hold true if \(U\in\cS\) has a multiple point 
\(p\) on the boundary \(\partial D\). Consider the harmonic function
\[
\Psi(x_1,x_2)=\left(x_2-(x_1+1)\right)
\left(x_2+(2+\sqrt{3})(x_1+1)\right)\left(x_2+(2-\sqrt{3})(x_1+1)\right)
\]
which nodal lines are represented in the following figure.
\begin{figure}[h!]
\begin{center}
\begin{tikzpicture}
\begin{axis}
[axis equal,axis lines=middle,enlargelimits,
xtick=\empty,
ytick=\empty]
\filldraw(0,100) circle (3pt);
\draw(187,50) circle (3pt);
\draw (14,49) circle (3pt);
\draw (100,200) circle (3pt);
\put(15,90){\(p_1\)};
\put(155,30){\(p_3\)};
\put(32,35){\(p_2\)};
\put(100,160){\(p_4\)};
\put(20,60){\(\gamma_1\)};
\put(155,130){\(\gamma_3\)};
\put(70,10){\(\gamma_2\)};
\put(38,130){\(\gamma_4\)};
\addplot
[domain=-1:1,samples=40,smooth,variable=\t,
 dotted]({-1},{t});
\addplot
[domain=-1:0,samples=40,smooth,
 thick]{x+1};
\addplot
[domain=-1:-0.86,samples=40,smooth,
 thick]{-(2+sqrt(3))*(x+1)};
\addplot
[domain=-1:0.86,samples=40,smooth,
 thick]{(-2+sqrt(3))*(x+1)};\addplot
[domain=0:360,variable=\t,
 samples=40,smooth,thick]
({cos(t)},{sin(t)});
\end{axis}
\end{tikzpicture}
\caption{One 4-point on the boundary.}\label{quadruplofrontiera}\label{fig0}
\end{center}
\end{figure}
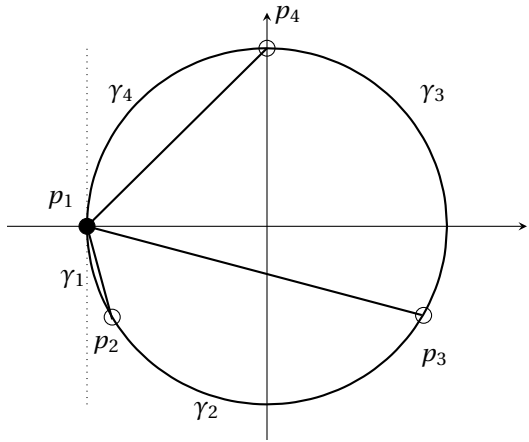 \\
On the boundary of the unit disk \(\Psi\) vanishes in \(p_1=(-1,0)\), 
\(p_{2,3}=(\mp\sqrt{3}/2,-1/2)\), \(p_4=(0,1)\). 
Denote by \(\gamma_i\) the arcs of positivity of \(\Psi|_{\partial D}\) 
and by \(\omega_i\) the sets of positivity of \(\Psi\) between 
\(\gamma_i\) and the lines where \(\Psi\) vanishes.
If we assume  \(\phi_i=\Psi\) on \(\gamma_i\), \(i=1,...,4\), 
then \(\Phi=(\phi_1,\phi_2,\phi_2,\phi_4)\)  is an admissible 
datum. Defining \(u_i=|\Psi|\) in \(\om_i\), then 
\(U=(u_1,u_2,u_3,u_4)\) has in \(p_1\in\partial D\) a 4-point, 
that is a point with multiplicity 4, and it can be easily seen 
that the angles between the interfaces are different 
(see also figure \ref{quadruplofrontiera}). \\
On the contrary the function \(U(x_1,x_2)=|x_1x_2|\) belongs to 
\(\cS\) and has at the origin an interior point with multiplicity 
4 with orthogonal nodal lines. \\
Note that in paper \cite{tvz} it is considered a segregation model in 
which the competition between the species is not symmetric, in that case 
the authors are able to prove that the multiple points cannot occur 
on the boundary of the domain.
\eoss

\section{Main results and proofs}

Define the set of points of multiplicity greater than or equal to 
$h\in\enne$
\[
\cZ_h(U)=\{x\in \overline D: m(x)\geq h\}\,.
\]
The set $\cZ_h(U)$ consists of isolated points \cite[Theorem 9.13]{ctv3}.

\begin{proposition}\label{P3.1}
Let \(U\in\cS\), then only one of the following statement is satisfied \\
{\rm i.} $\cZ_3(U)$ consists of one point $a_U\in \overline D$ 
	such that $m(a_U)=4$, \\
{\rm ii.} $\cZ_3(U)$ consists of two points $a_U, b_U\in \overline D$,
	$a_U\neq b_U$, with $m(a_U)=m(b_U)=3$.
\end{proposition}

\dimo The set \(\cZ_3(U)\) is nonempty. Indeed, if \(\cZ_3(U)=\emptyset\) 
then the interfaces between any two densities do not intersect. 
Since \(\phi\) vanishes in  exactly \(4\) points, there are only two 
interfaces which are nonempty, with endpoints on the boundary and do not 
intersect. This contradicts the fact that \(U\) defines four nodal regions.\\
Let \(a_U\in\overline{B}\) with \(m(a_U)=4\). Then any neighborhood 
of $a_U$ contains points of every $\omega_i$, and every non empty 
\(\Gamma_{ij}\) is connected, starts from \(\pa B\) and satisfies 
\(\overline{\Gamma}_{ij}\ni a_U\). Suppose that there exists 
\(b_U\neq a_U\) with $m(b_U)\geq 3$. Then \(b_U\) belongs to an interface 
between two densities, say \(b_U\in\Gamma_{12}\). It follows that there 
exists a neighborhood \(I\) of \(b_U\) such that \(I\cap \om_h\neq \emptyset\), 
\(h=1,2\), and \(I\) has positive distance from \(\om_k\), \(k=3,4\).  
On the other hand, since \(m(b_U)\geq 3\), \(I\cap \om_h\neq \emptyset\) 
for three different values of \(h\). This contradicts the fact that every 
$\om_i$ is connected.\\
If there aren't points with multiplicity $4$, then there are at least 
two points \(a_U\neq b_U\) of multiplicity $3$. Suppose that there 
exists \(c_U\) with \(m(c_U)=3\), \(c_U\neq a_U, b_U\). Then \(c_U\) 
belongs to  an interface between two densities, for example 
\(c_U\in\Gamma_{12}\).  Since \(m(c_U)=3\) there exists a neighborhood 
\(I\) of \(c_U\) such that \(I\cap \om_h\neq \emptyset\) for three 
different values of \(h\), for example \(h=1,2,3\).
Then \(\om_3\) is not connected, a contradiction. \finedim

\boss
All the 4-partitions of the disk \(D\) generated by the admissible 
boundary datum \(\Phi=(\phi_1,\phi_2,\phi_3,\phi_4)\) are topologically 
equivalent to one of the two configurations in figure \ref{figtype}
where \tikz\draw(0,0)circle(3pt); stands for an isolated zero of 
the boundary datum \(\Phi\), \tikz\filldraw[gray](0,0)circle(2pt); 
is a 3-point and \tikz\filldraw[black](0,0)circle(3pt); indicates 
a 4-point.
\begin{figure}[h!]
\begin{center}
\begin{tikzpicture}
\draw[thick] (0,0) circle (2.5cm);
\draw (-2.5,0) -- (2.5,0);
\draw (0,-2.5) -- (0,2.5);
\draw (-2.5,0) circle (3pt);
\draw (2.5,0) circle (3pt);
\draw (0,-2.5) circle (3pt);
\draw (0,2.5) circle (3pt);
\filldraw[black] (0,0) circle (3pt);
\end{tikzpicture}
\qquad
\begin{tikzpicture}
\draw[thin] (0,0) circle (2.5cm);
\draw (-2.5,0) -- (-0.7,-0.7);
\draw (0,-2.5) -- (-0.7,-0.7);
\draw (2.5,0) -- (0.7,0.7);
\draw (0,2.5) -- (0.7,0.7);
\draw (-0.7,-0.7) -- (0.7,0.7);
\draw (-2.5,0) circle (3pt);
\draw (2.5,0) circle (3pt);
\draw (0,-2.5) circle (3pt);
\draw (0,2.5) circle (3pt);
\filldraw[gray] (-0.7,-0.7) circle (2pt);
\filldraw[gray] (0.7,0.7) circle (2pt);
\end{tikzpicture}
\end{center}
\caption{Configurations with one 4-point (on the left) and two 
3--points (on the right).}
\label{figtype}
\end{figure}
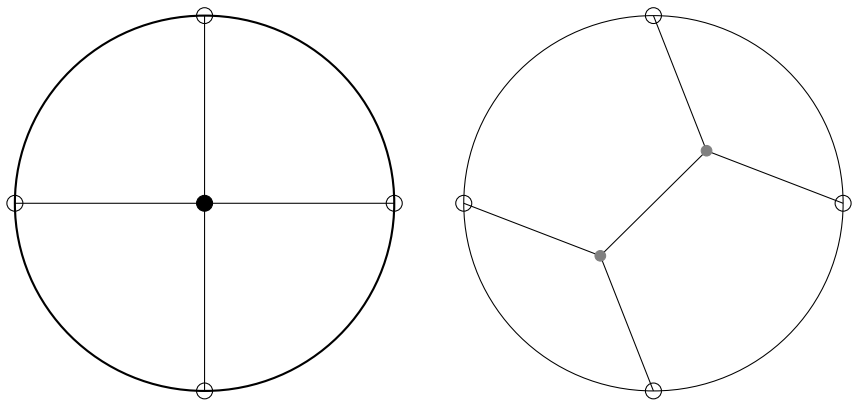 \\
In particular the multiple points can also lie on the 
boundary of \(D\). In figures \ref{fig0} and \ref{fig1} we show, 
respectively, configurations with a \(4\)-point on \(\partial D\) 
and with two \(3\)-points on \(\partial D\). The multiple points 
of \(U\) lie on \(\ov{D}\), and it is possible to obtain either 
one 4-point or two 3-points.
\eoss

Consider the boundary value problem
\beq\label{armonico}
\begin{array}\{{cl}.
-\Delta\psi= 0 & \hbox{ in }D \\
	\psi=\phi & \hbox{ on }\pa D \\
\end{array}
\eeq

\bteo\label{aleharmonic}
Let \(\Phi=(\phi_1,\phi_2,\phi_3,\phi_4)\) be an admissible datum.
The harmonic function \(\psi_a\) which solves \eqref{armonico} with 
boundary datum \(\phi^a= \sum_{j=1}^4 (-1)^j\phi_j\) possesses at most 
one critical point \(p\) in \(D\) such that \(\psi_a(p)=0\).
\eteo

\dimo Since \(\Phi\) is an admissible datum, the solution of 
\eqref{armonico} with boundary datum \(\phi^a\) vanishes at exactly 
four points on \(\partial D\), each arc 
\(\Gamma_j=\{\phi_j>0\}\subset\partial D\) is connected and 
\(\psi_a\) has different signs on adjacent arcs.

\begin{figure}[h!]
\begin{center}
\begin{tikzpicture}
\draw[thin] (0,0) circle (1.5cm);
\draw (-1.5,0) arc (270:359:1.5);
\draw (0,-1.5) arc (180:90:1.5);
\draw (-1.5,0) circle (2pt);
\draw (1.5,0) circle (2pt);
\draw (0,-1.5) circle (2pt);
\draw (0,1.5) circle (2pt);
\put(-40,-58){\(a)\)};
\put(1,1){\(\psi_a>0\)};
\put(35,35){\(\phi_4\)};
\put(-45,35){\(-\phi_1\)};
\put(-45,-35){\(\phi_2\)};
\put(30,-35){\(-\phi_3\)};
\end{tikzpicture}
\qquad
\begin{tikzpicture}
\draw[thin] (0,0) circle (1.5cm);
\draw (-1.5,0) -- (1.5,0);
\draw (0,-1.5) -- (0,1.5);
\draw (-1.5,0) circle (2pt);
\draw (1.5,0) circle (2pt);
\draw (0,-1.5) circle (2pt);
\draw (0,1.5) circle (2pt);
\filldraw[black] (0,0) circle (2pt);
\put(4,15){\(\psi_a>0\)};
\put(-30,-20){\(\psi_a>0\)};
\put(-40,-58){\(b)\)};
\put(35,35){\(\phi_4\)};
\put(-45,35){\(-\phi_1\)};
\put(-45,-35){\(\phi_2\)};
\put(30,-35){\(-\phi_3\)};
\end{tikzpicture}
\qquad
\begin{tikzpicture}
\draw[thin] (0,0) circle (1.5cm);
\draw (-1.5,0) -- (1.5,0);
\draw (-1.5,0) -- (0,1.5);
\draw (-1.5,0) -- (0,-1.5);
\filldraw[black] (-1.5,0) circle (2pt);
\draw (1.5,0) circle (2pt);
\draw (0,-1.5) circle (2pt);
\draw (0,1.5) circle (2pt);
\put(-40,-58){\(c)\)};
\put(4,15){\(\psi_a>0\)};
\put(35,35){\(\phi_4\)};
\put(-45,35){\(-\phi_1\)};
\put(-45,-35){\(\phi_2\)};
\put(30,-35){\(-\phi_3\)};
\end{tikzpicture}
\end{center}
\caption{The level  set \(\Gamma=\{x\in D:\psi_a(x)=0\}\) and the 
nodal components of \(\psi_a\), where the function \(\psi_a\) 
solves \eqref{armonico} with boundary datum 
\(\phi^a= \sum_{j=1}^4 (-1)^j\phi_j\). Three different situations can occur. }\label{figtypebis}
\end{figure}
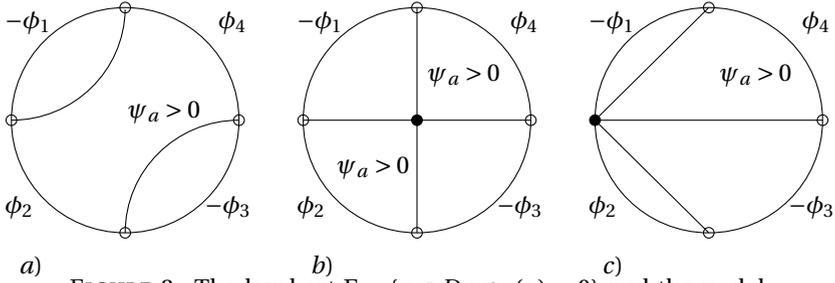

\begin{figure}[h]
\begin{center}
\begin{tikzpicture}
\draw[thick] (0,0) circle (2.5cm);
\draw (-2.5,0) -- (0,2.5);
\draw (0,-2.5) -- (2.29,1);
\draw (-2.5,0) -- (0,-2.5);
\draw (0,2.5) circle (3pt);
\draw (2.29,1) circle (3pt);
\put(60,60){\(\phi_4\)};
\put(-60,60){\(\phi_1\)};
\put(-60,-60){\(\phi_2\)};
\put(60,-60){\(\phi_3\)};
\put(-84,0){\(p\)};
\put(0,-84){\(q\)};
\filldraw[gray] (0,-2.5) circle (2pt);
\filldraw[gray] (-2.5,0) circle (2pt);
\put(-90,-90){\(a)\)};
\end{tikzpicture}
\qquad
\begin{tikzpicture}
\draw[thick] (0,0) circle (2.5cm);
\draw (-2.5,0) -- (0,2.5);
\draw (0,-2.5) -- (2.5,0);
\draw (-2.5,0) -- (2.5,0);
\draw (0,2.5,0) circle (3pt);
\draw (0,-2.5) circle (3pt);
\put(60,60){\(\phi_4\)};
\put(-60,60){\(\phi_1\)};
\put(-60,-60){\(\phi_2\)};
\put(60,-60){\(\phi_3\)};
\put(-84,0){\(p\)};
\put(75,0){\(q\)};
\filldraw[gray] (2.5,0) circle (2pt);
\filldraw[gray] (-2.5,0) circle (2pt);
\put(90,-90){\(b)\)};
\end{tikzpicture}
\end{center}
\caption{Configuration with 3-points on \(\partial D\).}
\label{fig1}
\end{figure}
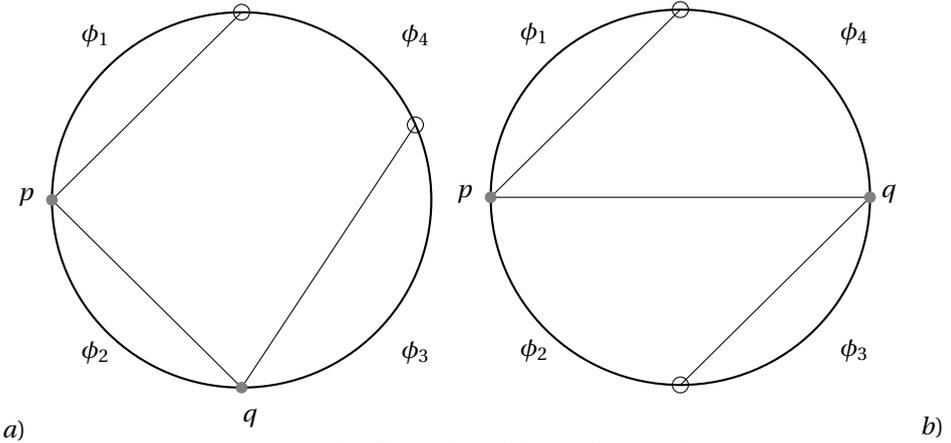

Since a harmonic function does not admit closed level lines, the set 
\(\Gamma=\{x\in D:\psi_a(x)=0\}\) has no closed loop. We infer that 
\(\psi_a\) has alternate positive or negative sign on \(3\) or \(4\) 
adjacent sets: the nodal components of \(\psi_a\) (see figure \ref{figtypebis}).

By standard theory of harmonic functions the zero set of \(\psi_a\) 
around a critical point at level 0 is made by (at least) \(4\) half-lines, 
meeting with equal angles. We infer that locally around each critical 
point at level \(0\) the function \(\psi_a\) defines \(4\) nodal
components. Such components are exactly the \(4\) nodal regions of 
\(\psi_a\), since we have already pointed out that \(\Gamma\) does
not contain closed loop.

As a consequence, if the function \(\psi_a\) has \(3\) nodal components 
(figure \ref{figtypebis}a), there are no critical points at level \(0\).  
If the nodal components of \(\psi_a\) are \(4\), two different situations 
can occur.\\
1) The set \(\Gamma\) is made by \(2\) simple arcs with endpoints on the 
boundary, intersecting in a point \(p\in D\) (figure \ref{figtypebis}b). 
Let \(q\) be a critical point  at level \(0\), with \(q\neq p\).  
Then one can construct a closed Jordan curve passing through \(p\) 
and \(q\), belonging to the positivity (or negativity) set of \(\psi_a\),  
which disconnects \(\partial D\) from a negative (or positive) nodal 
components. This is impossible because  \(\Gamma\) has no closed loop.
\\
2) The set \(\Gamma\) is made by \(4\) simple arcs with endpoints on the 
boundary, intersecting in \(p\in \partial D\) (figure \ref{figtypebis}c). 
If there exists a critical point \(q\) at level \(0\), with \(q\neq p\) 
with the same arguments used in 1)  we get a contradiction.

\begin{remark}\rm
The existence of a critical point in theorem \ref{aleharmonic} is not 
guaranteed. Consider, for example, the harmonic polynomial
\[
\psi_a(x_1,x_2)=4(x_1^2-x_2^2)+10 x_1+7
\]
On the boundary of the unit disk
\[
\psi_a\big|_{\partial D}=\Phi^a(\theta)=(4 \cos(\theta)+3)(2\cos(\theta)+1)
\]
Denoting by \(\theta_0\in (\pi/2,\pi): \cos(\theta_0)=-3/4\), we have
\[
\begin{aligned}
\Phi^a(\theta)&<0 \quad {\rm in}\quad(2\pi/3,\theta_0)\cup 
  (-\theta_0,-2\pi/3 )\\
\Phi^a(\theta)&>0 \quad {\rm in}\quad( -2\pi/3,2\pi/3)\cup 
  (\theta_0,-\theta_0)
\end{aligned}
\]
Hence \(\Phi=|\Phi^a|\) is admissible and the unique critical point 
of \(\psi_a\) is \((-5/4,0)\) which does not belong to \(D\).
Conditions which ensure the existence of critical points for solutions 
of two dimensional elliptic equation are given in  \cite{w} and \cite{a}.
\end{remark}

\boss
Recall that all the critical points of \(\psi_a\) in \(D\) are always 
saddle points, since harmonic functions in \(D\) assume local maximum 
or local minimum only on the boundary of the ball. 
If \(x_0\) is an internal critical point then, in a neighborhood of \(x_0\), 
the level line \(\left\{x\in D: \psi_a(x)=\psi_a(x^0)\right\}\) is made 
by \(2\) simple arcs intersecting in \(x_0\) (see \cite[Remark 1.2]{a}). 
\eoss

\begin{proposition}\label{iffharmonic} 
Let \(\Phi\) be an admissible datum. If the set \(\mathcal{S}\) contains 
a function \(U\) with a 4-point \(a_U\) in \(\overline{D}\) then \(U=|\psi_a|\), 
where \(\psi_a\) is the harmonic function such that 
\(\psi_a=\phi^a= \sum_{j=1}^4 (-1)^j\phi_j\) on \(\partial D\).
\end{proposition}

\dimo Let \(U\) be an element of  \(\mathcal{S}\)  with a 4-point \(a_U\), 
then there exist three or four arcs connecting \(a_U\) to any of the 
isolated zeros of the boundary datum (note that there are three arcs 
if and only if \(a_U\in\partial D\)). Then the function 
\(\psi_a=\sum_{j=1}^4(-1)^ju_j\) is harmonic in \(D\setminus\{a_U\}\) 
(see \cite[Proposition 6.3]{ctv3}), moreover \(\psi_a\) is bounded so,
by Schwarz's removable singularity principle (see \cite[Proposition 11.1]{p}) 
\(\psi_a\) is harmonic in \(D\) and, by construction, \(U=|\psi_a|\) in \(D\) 
and \(U=|\Phi^a|=\Phi\) on \(\partial D\).\finedim

\begin{proposition}\label{3.9}
\(\mathcal{S}\) contains only solutions of one type. In other words it is 
not possible that there exist \(U,W\in\mathcal{S}\) such that \(U\) 
possesses one 4-point \(a_U\) and \(W\) possesses two 3-points \(a_W,b_W\).
\end{proposition}

\dimo This statement is an easy consequence of the uniqueness results 
contained in \cite{ab,wz}, however we present here a different proof. \\
Suppose there exist \(U,W\in\mathcal{S}\) such that \(U=(u_1,...,u_4)\) 
possesses one 4-point \(a_U\) and \(W=(w_1,...,w_4)\) possesses two 3-points 
\(a_W,b_W\) for an admissible boundary datum \(\Phi=(\phi_1,...,\phi_4)\). 
Since \(W\) has two 3-points we can assume that there are two nodal set having 
positive distance, without loss of generality we can think that 
\(\overline{\om}_{W,1}\cap\overline{\om}_{W,3}=\emptyset\), where 
\(\om_{W,1}=\{x\in D: w_1(x)>0\}\) and \(\om_{W,3}=\{x\in D: w_3>0\}\)
(see, for example, figure \ref{figtype}). Consider the following functions
\[
w^*=w_1-w_2+w_3-w_4 \,,\qquad  \psi_a=u_1-u_2+u_3-u_4\,.
\]
By proposition \ref{iffharmonic} we know that \(\psi_a\) is an harmonic 
function in \(D\) such that \(\psi_a=\phi^a\) on \(\partial D\), moreover
\(w^*\) is a superharmonic function, since is a gluing of harmonic functions 
which become superharmonic on the arc linking the \(a_W\) and \(b_W\) (see 
\cite[Proposition 6.3]{ctv3}, and also \cite[Chapter 12]{p}). \\
Then the function \(z(x_1,x_2)=w^*(x_1,x_2)-\psi_a(x_1,x_2)\), with 
\((x_1,x_2)\in D\), is a superharmonic function such that \(z(x_1,x_2)=0\) 
if \((x_1,x_2)\in\partial D\), and, by maximum principle, it follows 
\(z(x_1,x_2)\geq0\), that is
\begin{equation}\label{ineq}
w^*(x_1,x_2)\geq\psi_a(x_1,x_2) \qquad \forall(x_1,x_2)\in D\,.
\end{equation}
Inequality \eqref{ineq} implies that the positivity set of \(\psi_a\) is
contained in the positivity set of \(w^*\). Precisely it follows that
\(\om_{U,1}\subseteq\om_{W,1}\) and \(\om_{U,3}\subseteq\om_{W,3}\)
and the contradiction is reached since 
\(\overline{\om}_{U,1}\cap\overline{\om}_{U,3}=a_U\)
but \(\overline{\om}_{W,1}\cap\overline{\om}_{W,3}=\emptyset\).
\finedim

\begin{proposition}\label{3.10}
Let \(U,W\in\mathcal{S}\) be functions with only one 4-point 
\(a_U,a_W\in\ov{D}\), then \(U=W\).
\end{proposition}

\dimo Assume there exist two different \(U,W\in\mathcal{S}\) with a 4-point, then 
there exits three or four arcs connecting \(a_U\) to any of the isolated 
zeros of the boundary datum (there are three arcs if and only if 
\(a_U\in\partial D\)), and the same is true for \(a_W\).
In any case we can consider the functions \(u_a=\sum_{j=1}^4(-1)^j u_j\) 
and \(w_a=\sum_{j=1}^4(-1)^j w_j\), by proposition 6.3 in \cite{ctv3}, both 
functions are solution to \eqref{armonico} with boundary datum 
\(\phi^a=\sum_{j=1}^4(-1)^j\phi_j\), but it is well known that \eqref{armonico} 
possesses only one harmonic solution, that is \(u_a\equiv w_a\), 
and this means that \(U=W\).\finedim

Let \(U\in\cS\) be given with trace \(\Phi=(\phi_1,...,\phi_4)\) and 
3-points on \(\partial D\). 
The next proposition shows that, similarly to proposition \ref{iffharmonic},  
we can construct a harmonic function closely related to \(U\), which has 
different signs on adjacent nodal regions.

\begin{proposition}\label{3.9bis}
Suppose that \(U\in\cS\) is a configuration with two 3--points
on the boundary, that is there are \(p,q\in \partial D\) such that 
\(m(p)=m(q)=3\). Then we can construct a harmonic function \(\Xi_a\) 
such that \(U=|\Xi_a|\) in \(D\).
\end{proposition}

\dimo We denote by \(\Sigma_j=\{x\in\partial D: \phi_j(x)>0\}\), \(j=1,...,4\). 
Since \(\Phi\) is admissible we have \(\Sigma_i\cap \Sigma_j=\emptyset\), 
\(i\neq j\), and \(\partial D=\cup_{j=1}^4 \overline{\Sigma}_j\).  
Two different configurations are possible (see figure \ref{fig1}) and,  
in both cases, we construct a function which is alternatively positive 
and negative on the nodal region of \(U\) and is harmonic in \(D\).\\
a) There exists an index \(i\) such that \(p\) and \(q\) belong to 
\(\overline{\Sigma}_i\) that is \(\phi_i(p)=\phi_i(q)=0\), say \(i=2\). 
Then there exists  \(j\neq i\) such that \(p,q\not\in \overline{\Sigma}_j\), 
say \(j=4\) (see figure \ref{fig1}, on the left). 
Then the function \(\Xi_a=u_4-(u_1+u_2+u_3)\) has different signs on 
adjacent nodal regions. By \cite[Proposition 6.3]{ctv3}, \(\Xi_a\) is 
harmonic in \(D\setminus \Omega_{12}\), in \(D\setminus \Omega_{13}\) 
and in \(D\setminus \Omega_{23}\), where 
\(\Omega_{hk}=\overline{\om}_h\cup\overline{\om}_k\). 
As \(\Omega_{12}\cap\Omega_{23}\cap\Omega_{13}=\{p,q\}\) we deduce that 
\(\Xi_a\) is harmonic in \(D\). By construction, \(U=|\Xi_a|\) in \(D\).\\
b) \(p,q\) do not belong to the  same \(\overline{\Sigma}_i\) for 
any \(i\) that is they are not consecutive. Suppose, e.g., that 
\(p\in \overline{\Sigma}_1\), \(q\in \overline{\Sigma}_3\), (see figure 
\ref{fig1}, on the right).
Then the function \(\Xi_a=u_1+u_2-u_3-u_4\) has different signs on 
adjacent nodal regions and, by \cite[Proposition 6.3]{ctv3}, \(\Xi_a\) is 
harmonic in \(D\setminus \Omega_{23}\), in \(D\setminus \Omega_{13}\) and 
in \(D\setminus \Omega_{14}\). 
As \(\Omega_{23}\cap\Omega_{13}\cap\Omega_{14}=\{p,q\}\) we deduce that 
\(\Xi_a\) is harmonic in \(D\). The function  \(\Xi_a\) solves the problem 
\eqref{armonico} with boundary datum \(\phi_1+\phi_2- \phi_3-\phi_4\) 
and, by construction,  \(U=|\Xi_a|\) in \(D\).
\finedim

\begin{remark}\rm If \(U\in\cS\) produces a 3-points configuration 
with at least a \(3\)-point \(p\in D\) then proposition \ref{3.9} 
does not hold true. Indeed, since \(m(p)=3\), there is a neighborhood 
\(I\) around  \(p\) such that \(I\cap \om_h\neq \emptyset\) for three 
different values of \(h\), and three interfaces with endpoint \(p\) 
are not empty.  Then it is not possible to construct an harmonic 
function \(\Xi_a\) such that \(U=|\Xi_a|\) with different signs on 
adjacent nodal regions around \(p\).
\end{remark}

Let us remark that, if \(p\in D\) with \(m(p)=4\) then
\[
U(r,\th)=r^{2} \left| \cos\left( 2 (\th +\th_0)\right)
\right|+o(r^{2})
\]
as $r\to 0$, where $(r,\th)$ is a system of polar coordinates around 
\(p\). Hence \(U(p)=0\) and \(\nabla U(p)=0\). Proposition \ref{iffharmonic} 
can be inverted if \(p\in D\).

\begin{proposition}\label{4point}
Let \(\Phi\) be an admissible boundary datum and suppose that the 
related harmonic function \(\psi_a\), solution to \eqref{armonico} 
with boundary datum \(\phi^a\), has a critical point \(p\in D\) 
(a saddle) such that \(\psi_a(p)=0\), then \(p\) is a 4-point 
for the function \(U=|\psi_a|\in\cS\).
\end{proposition}

\dimo 
We introduce the transformation 
\begin{equation}\label{conforme}
x=T_p(\zeta)=\frac{\zeta+p}{\overline{p}\zeta+1}\,.
\end{equation}
Here we identificate the complex numbers \(x=x_1+i x_2\) and \(\zeta=\zeta_1+i\zeta_2\) with the 
points \((x_1,x_2)\) and \((\zeta_1,\zeta_2)\in\erre^2\), respectively. \(T_p\) is a conformal map which maps 
the unit disk \(\overline{D}\) into itself such that \(T_p(\partial D)=\partial D\) 
and \(T_p(0)=p\). Set
\begin{equation}\label{TR}
\widetilde \Psi(\zeta)=\Psi_a(T_p(\zeta)),\qquad
\widetilde \Phi^a(\zeta)=\Phi^a(T_p(\zeta))\,.
\end{equation}
Then \(\widetilde \Psi\) solves the problem
\beq\label{DP2}
\begin{array}\{{cl}.
-\Delta\widetilde \Psi= 0 & \hbox{ in }D\,, \\
	\widetilde \Psi=\widetilde \Phi^a & \hbox{ on }\pa D\,. \\
\end{array}
\eeq
Denote by
\(u(\z)=\re(T_p(\z))\) and \( v(\z)=\im(T_p(\z))\). By differentiation we get
\[
\nabla_\z  \widetilde \Psi (\z)= \left(\begin{matrix}
\partial_{\z_1} u(\z)& \partial_{\z_1} v(\z)\\
\partial_{\z_2} u(\z)& \partial_{\z_2} v(\z)
\end{matrix}\right)
\nabla_x   \Psi (T_p(\z))\,.
\]
From the relations
\[
\frac{\partial }{\partial \z_1}T_p(\z)=\frac{1-|p|^2}{(1+\overline p \z)^2} 
\qquad
\frac{\partial }{\partial \z_2}T_p(\z)=i\frac{1-|p|^2}{(1+\overline p \z)^2}\,.
\]
we obtain 
\[
\left(\begin{matrix}
\partial_{\z_1} u(\z)& \partial_{\z_1} v(\z)\\
\partial_{\z_2} u(\z)& \partial_{\z_2} v(\z)
\end{matrix}\right)\Bigg|_{\z=0}=\left(\begin{matrix}
(1-|p|^2)&0\\0&(1-|p|^2)
\end{matrix}\right)\,.
\]
Hence
\begin{equation}\label{nabla}
\nabla_\z  \widetilde \Psi (0)= \left(\begin{matrix}
(1-|p|^2)&0\\0&(1-|p|^2)
\end{matrix}\right)
\nabla_x\Psi_a (p)\,.
\end{equation}
We can write the Fourier expansion of \(\widetilde \Psi\)
\[
\widetilde \Psi(r,\theta)=\frac{A_0}{2}
+\sum_{k=1}^\infty (A_k \cos(k\theta)+B_k \sin(k\theta))r^k\,,
\qquad \zeta=(r,\theta)\,.
\]
The hypothesis 
\(
\widetilde \Psi(0)=\Psi(p)=0
\)
implies that \(A_0=0\) whereas the hypothesis  \(\nabla_x \Psi(p)=(0,0)\) 
gives  \(\nabla_\z \widetilde \Psi(0)=(0,0)\) and then  \(A_1=B_1=0\).
Moreover \((A_2,B_2)\neq (0,0)\)  because \(\widetilde \Psi\) has exactly 
four zeroes  on the boundary. Therefore 
\(U(x)=|\Psi_a(x)|=|\widetilde \Psi (T^{-1}_p(x))|\)
is nonnegative, satisfies the boundary datum \(\Phi^a\) and has 
exactly four nodal regions. This function generates an element of \(\cS\), with datum \(\Phi\) and quadruple point \(p\) (see also \cite[Lemma 3.2]{ctv4}).
\finedim

\begin{proposition}\label{CNeS}
Suppose that  \(\Phi=(\phi_1,\phi_2,\phi_3,\phi_4)\) is an admissible datum.
The solution \(U\) of the minimum problem \eqref{minimo} with 
admissible datum \(\Phi\) has a  4-point in \(p\in {D}\) if and only if
\begin{eqnarray}\label{C1}
&\displaystyle\sum_{j=1}^4 (-1)^j \int_{\partial D} \phi_j\left(
	\frac{\zeta+p}{\overline{p}\, \zeta+1}\right) ds_\zeta&=0\,,\\\label{C2}
&\displaystyle\sum_{j=1}^4 (-1)^j \int_{\partial D} \phi_j\left(
	\frac{\zeta+p}{\overline{p}\, \zeta+1}\right) \zeta_r ds_\zeta
	&=0\,, \qquad r=1,2\,.
\end{eqnarray}
Here we use the notation \(\zeta=\zeta_1+i \zeta_2\).
\end{proposition}
\dimo
Suppose that \(U\) has a  4-point in \(p\in {D}\). 
For proposition \ref{iffharmonic} we have \(U=|\Psi_a|\) where \(\Psi_a\) 
solves \eqref{armonico} with boundary datum \(\Phi^a=\sum_{j=1}^4 (-1)^j \phi_j\), and, by \eqref{local},
\(U(p)=0\), \(\nabla U(p)=(0,0)\). Hence \(\nabla \Psi_a (p)=(0,0)\).
We introduce the transformation \eqref{conforme} and define 
\(\widetilde \Psi\) and \(\widetilde \Phi^a\) according to \eqref{TR}.
Then \(\widetilde \Psi\) solves \eqref{DP2}
and, by hypothesis, the origin is a 4-point for 
\(\widetilde U(\zeta)=|\widetilde \Psi(\zeta)|\). By the Poisson 
integral formula 
\begin{equation}\label{poisson}
\widetilde \Psi(\zeta)=\frac{1-|\zeta|^2}{2\pi}
\int_{\partial D} \frac{\widetilde \Phi^a (\eta)}{|\zeta-\eta|^2}ds_\eta
\end{equation}
and it belongs to \(C^2(D)\cap C^0(\overline{D})\) (cf. \cite[(2.27)]{gt}).
We deduce that
\[
0=\widetilde \Psi(0)=\frac{1}{2\pi}\int_{\partial D} 
\frac{\widetilde \Phi^a (\eta)}{|\eta|^2}ds_\eta=
\frac{1}{2\pi}\int_{\partial D} \widetilde \Phi^a (\eta)ds_\eta=
\frac{1}{2\pi}\int_{\partial D}  \Phi^a (T_p(\eta))ds_\eta
\]
which gives \eqref{C1}. By \eqref{nabla} we have 
\(\nabla\widetilde{\Psi}(0)=(0,0)\).
By direct differentation of the Poisson integral
\[
\frac{\partial}{\partial \zeta_r} \widetilde\Psi (0)= \frac{1}{\pi} \int_{\partial D} \widetilde \Phi^a(\eta)\eta_r ds_\eta =
 \frac{1}{\pi} \int_{\partial D}  \Phi^a(T_p(\eta))\eta_r ds_\eta
=0,\quad r=1,2
\]
which gives \eqref{C2}.

Conversely, suppose that \eqref{C1} and \eqref{C2} are valid. 
Let \(\Psi_a\) be the solution of \eqref{armonico} with boundary 
datum \(\Phi^a=\sum_{j=1}^4 (-1)^j \phi_j\). 
Then the function \(\widetilde \Psi(\zeta)=\Psi_a(T_p(\zeta))\) solves 
\eqref{DP2} and is given by \eqref{poisson}.
By \eqref{C1}  we have 
\[
\begin{aligned}
&\widetilde \Psi(0)=
\frac{1}{2\pi}\int_{\partial D}  \Phi^a (T_p(\eta))ds_\eta=0\,.
\end{aligned}
\]
Hence  \(\Psi_a(p)=\Psi_a(T_p(0))=\widetilde \Psi(0)=0\). \\
Moreover, by \eqref{C2},
\[
\frac{\partial}{\partial \zeta_s} \widetilde\Psi (0)= \frac{1}{\pi} \int_{\partial D} \widetilde \Phi^a(\eta)\eta_s ds_\eta =
 \frac{1}{\pi} \int_{\partial D}  \Phi^a(T_p(\eta))\eta_r ds_\eta=0,\quad r=1,2.
\]
From \eqref{nabla} we get
\[
\nabla_x \Psi_a(p)= \left(\begin{matrix}
(1-|p|^2)^{-1}&0\\0&(1-|p|^2)^{-1}
\end{matrix}\right)
\nabla_\zeta  \widetilde \Psi (0)=(0,0)\,.
\]
For proposition \ref{4point} we conclude that the function \(U=|\Psi_a|\) posses a 
4-point and belongs to \(\cS\).
\finedim

\boss\label{codim}
Conditions \eqref{C1},\eqref{C2} imply that an admissible datum generates a 4-point configuration if and only if \(3\) integral conditions are satisfied.
These conditions identify three hyperplanes in the space \(W^{1,\infty}(\partial D)\), then the boundary data which generate a 4-point configuration 
form (at most) a codimension \(3\) set in the space of admissible data. 
\eoss

If \(p\) is the origin, proposition \ref{CNeS} can be formulated as 
follows

\begin{corollary}\label{CNES2}
Let \(\Phi=(\phi_1,\phi_2,\phi_3,\phi_4)\) be an admissible boundary 
datum. 
Then \(U=|\psi_a|\), where \(\psi_a\) is the solution of 
\eqref{armonico} with boundary datum 
\(\phi^a\), has a 4-point in the origin  if and only if
\begin{equation}\label{cond12}
\int_{\pa D} \phi^a(y)\, ds_y=0
\qquad \hbox{and}\qquad
\int_{\pa D}y_j \phi^a(y)\, ds_y=0\quad j=1,2
\end{equation}
where  \(\phi^a=\sum_{j=1}^4 (-1)^j \phi_j\). 

\end{corollary}

\boss
It is easy to verify that conditions \eqref{cond12} hold if the admissible 
datum  \(\Phi=(\phi_1,\phi_2,\phi_3,\phi_4)\) satisfies the following 
symmetry conditions
\[
\phi_2(x_1,x_2)=\phi_1(-x_1,x_2) \quad
\phi_3(x_1,x_2)=\phi_1(-x_1,-x_2) \quad
\phi_4(x_1,x_2)=\phi_1(x_1,-x_2)
\]
for any \(x=(x_1,x_2)\in\partial D\).
\eoss

Proposition \ref{CNeS} gives necessary 
and sufficient conditions on the datum such that the solution of the 
minimum problem \eqref{minimo} generates a 4-point configuration with 
the quadruple point \(p\in D\) and  suggests that the most probable configurations 
in nature are the one with 3--points. However a 4-point 
configuration can be obtained as limit of 3--points configurations (see Remark \ref{codim}).
This is a consequence of the continuously dependence of the minimum in 
\eqref{minimo} on the boundary data (cf. \cite[Theorem 3.2]{ctv3}).
In fact we expect that, generically, the limit configuration
has two 3--points, with respect to a suitable
measure in the space of the admissible boundary data.

\begin{proposition}
Let \(U=(u_1,u_2,u_3,u_4)\) be the solution of the minimum problem 
\eqref{minimo} with admissible datum \(\Phi=(\phi_1,\phi_2,\phi_3,\phi_4)\), 
which generates a configuration with the 4-point \(q\in D\), then
there exists \(\eps_0>0\) and a continuum of functions
\(U^{(\eps)}=(u_1^{(\eps)},u_2^{(\eps)},\) \(u_3^{(\eps)},u_4^{(\eps)})\)
belonging to a class of segregated states, with admissible boundary datum,
such that \\
i. \(U^{(\eps)}\) is a two 3--points configuration, for any \(\eps\in(0,\eps_0)\),\\
ii. \(u_j^{(\eps)}\to  u_j\) in \(H^1(D)\), \(j=1,...,4\), as \(\eps\to 0\). 
\end{proposition}

\dimo   
Let \(\widetilde \Phi\) be an admissible datum which vanishes on the 
boundary \(\partial D\) in the same points of \(\Phi\)  and  such that 
the solution of the minimum problem \eqref{minimo} with boundary 
condition \(\widetilde \Phi\)  generates a configuration with  3-points 
\(t_1,t_2 \in D\). One among \(\Phi\pm \eps \widetilde \Phi\) is admissible 
for any \(\eps>0\). Let us assume that 
\(\Phi^{(\eps)}=\Phi+\eps \widetilde \Phi\) is admissibile.
Let \(\psi,\widetilde \psi\) and \(\psi^{(\eps)}\) be the solutions 
of \eqref{armonico} with boundary data \(\Phi^a, \widetilde \Phi^a\) 
and \(\Phi^a+\eps \widetilde \Phi^a\), respectively. It is obvious 
that \(\psi^{(\eps)}=\psi+\eps\widetilde \psi\).\\
If the solution \(U^{(\eps)}\) of \eqref{minimo}  with boundary datum 
\(\Phi^{(\eps)}\) has a 4-point \(Q\in D\) then, for proposition 
\ref{iffharmonic}, \(U^{(\eps)}=|\psi^{(\eps)}|\) and
\begin{equation}\label{condQ}
\psi^{(\eps)}(Q)=0,\qquad \nabla \psi^{(\eps)}(Q)= (0,0)\,.
\end{equation}
If \(Q=q\) then \(\psi(Q)=0,\nabla\psi (Q)=(0,0)\) but  at least one 
among \(\widetilde\psi(Q)=0\) and \(\nabla\widetilde\psi(Q)=(0,0)\) is not 
valid (because \(Q\) is not a 4-point for \(\widetilde U=|\widetilde\psi|\)). 
Hence \eqref{condQ} doesn't hold and we get a contradiction. 
It follows that \(Q\neq q\). \\
If \(\psi(Q)\neq 0\) and \(\widetilde \psi(Q)=0\) then  
\(\psi^{(\eps)}(Q)=\psi(Q)\neq 0\). Hence \(Q\) cannot be a 4-point 
for \(U^{(\eps)}\). If \(\psi(Q)= 0\) then \(\nabla \psi(Q)\neq (0,0)\) 
(because \(Q\) is not a 4-point for \(|\psi|\)). Since 
\(\nabla \psi^{(\eps)}(Q)=\nabla\psi(Q)+\eps \nabla \widetilde \psi(Q)\)  
there exists \(\eps_0>0\) such that  \(\nabla \psi^{(\eps)}(Q)\neq (0,0)\) 
for any  \(\eps<\eps_0\). It follows that \(U^{(\eps)}\) cannot generate 
a configuration with a 4-point if \(\eps<\eps_0\). 
For proposition \ref{P3.1} \(U^{(\eps)}\)  has two 3-points. \\
We have \(\Phi^{(\eps)}-\Phi=\eps \widetilde \Phi\to 0\) in 
\(H^{1/2}(\partial D)\). For the  continuous dependence of \(U^{(\eps)}\)  
on the boundary data (\cite[Theorem 3.2]{ctv3}) the result follows.
\finedim

\section{Some numerical examples}

\begin{figure}
\begin{center}
\includegraphics[scale=.28]{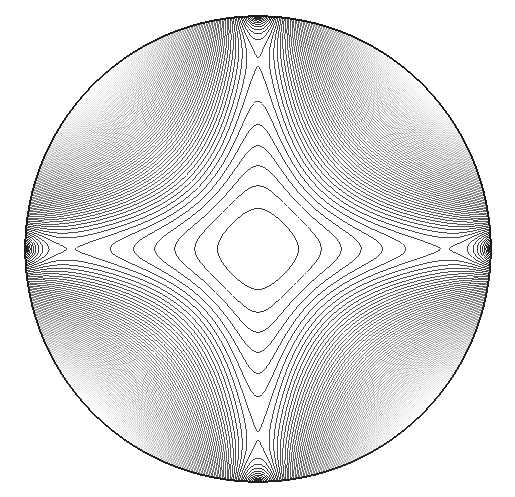}
\quad
\includegraphics[scale=.28]{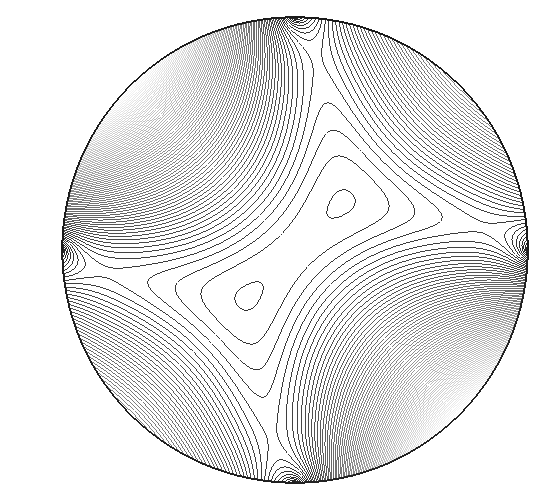}
\end{center}
\caption{Level lines of the approximate solutions of \eqref{diffeqk} with \(4\) differential equations and for large \(\mu\), showing a 4-point configuration (on the left) and a two 3-points configuration (on the right).}
\label{fig4}
\end{figure}

In \cite{ctv2} it is proved that the unique solution (see \cite{wz}) of the 
elliptic systems \eqref{diffeqk} converges to the solution of the segregation 
problem as \(\mu\to+\infty\) and, at least for \(k=2\) species, the rate of 
the convergence is \(\mu^{-1/6}\). Using this result we provide numerical 
simulations of solutions of the system 
\eqref{diffeqk} with \(4\) differential equations and for large \(\mu\).
Thanks to the FreeFEM++ software (a simple and powerful partial differential 
equation numerical solver, see \cite{ff}) we can obtain nodal partitions showing 
a 4-point solution (figure \ref{fig4}, on the left) and a two 3-points solution 
(figure \ref{fig4}, on the right).
\\
The pictures \ref{fig4}, reflecting the model configurations in figure \ref{figtype}, 
show some level lines of the approximate solutions of \eqref{diffeqk} with \(4\) equations and for large \(\mu\).
The function on the left side is the solution obtained with the
boundary data \(\phi_i(x_1,x_2)=15|x_1 x_2|\) for \(i=1,...,4\) and \(\mu=100\),
the picture on the right is obtained choosing \(\phi_i(x_1,x_2)=7|x_1 x_2|\) 
for \(i=1,3\), \(\phi_i(x_1,x_2)=15|x_1 x_2|\) for \(i=2,4\) and \(\mu=100\).
For the interested reader  we include the FreeFEM++ script which
produces the two 3-points configuration above. \\
{\small\begin{verbatim}
real R=1;
border gamma1(t=0,pi/2){x=R*cos(t);y=R*sin(t);}; 
border gamma2(t=pi/2,pi){x=R*cos(t);y=R*sin(t);}; 
border gamma3(t=pi,3*pi/2){x=R*cos(t);y=R*sin(t);}; 
border gamma4(t=3*pi/2,2*pi){x=R*cos(t);y=R*sin(t);}; 

mesh disk=buildmesh(gamma1(50)+gamma2(50)+gamma3(50)+gamma4(50));
plot(disk,wait=1); 

fespace H(disk,P2); 
H u1,u2,u3,u4,v1,v2,v3,v4,uold1=0,uold2=0,uold3=0,uold4=0,uT; 
real mu=100; 
real al=7,bt=15,ga=7,dt=15; 
func f1=al*abs(x*y); 
func f2=bt*abs(x*y); 
func f3=ga*abs(x*y); 
func f4=dt*abs(x*y); 

problem segregazione(u1,u2,u3,u4,v1,v2,v3,v4)= 
int2d(disk)(dx(u1)*dx(v1)+dy(u1)*dy(v1)) 
+int2d(disk)(dx(u2)*dx(v2)+dy(u2)*dy(v2)) 
+int2d(disk)(dx(u3)*dx(v3)+dy(u3)*dy(v3)) 
+int2d(disk)(dx(u4)*dx(v4)+dy(u4)*dy(v4))  
+int2d(disk)(mu*u1*(uold2+uold3+uold4)*v1) 
+int2d(disk)(mu*u2*(uold1+uold3+uold4)*v2) 
+int2d(disk)(mu*u3*(uold2+uold1+uold4)*v3) 
+int2d(disk)(mu*u4*(uold2+uold3+uold1)*v4)
+on(gamma1,u1=f1)+on(gamma2,u1=0)+on(gamma3,u1=0)+on(gamma4,u1=0) 
+on(gamma1,u2=0)+on(gamma2,u2=f2)+on(gamma3,u2=0)+on(gamma4,u4=0) 
+on(gamma1,u3=0)+on(gamma2,u3=0)+on(gamma3,u3=f3)+on(gamma4,u3=0) 
+on(gamma1,u4=0)+on(gamma2,u4=0)+on(gamma3,u4=0)+on(gamma4,u4=f4); 

real N=20; 
for(real t=1;t<N;t+=1){ 
segregazione; 
uold1=u1; 
uold2=u2; 
uold3=u3; 
uold4=u4; 
cout <<t<<endl;} 

uT=u1+u2+u3+u4; 
plot (uT,value=true,wait=true,fill=false,nbiso=80,grey=true);
\end{verbatim}}

\section{Conclusions}

The steady state of \(k\) competing species coexisting in the same area is 
governed by the competing-diffusion system \eqref{diffeqk}. 
Large interaction induces the spatial segregation of the species in the limit 
configuration as \(\mu\) tends to \(+\infty\), that is different densities have 
disjoint support.  We analysed the behaviour of the solutions in the case of \(4\) competing species as the competition 
rate tends to infinity and we proved that only two configurations are possible:  
or the interfaces between components meet in a 4-point that is a point where 
four components concur, or the interfaces between components meet in two 3-points, 
that is points where only three components concur (see figure \ref{fig4}).   

We characterized, for a given datum, the possible 4-point configuration by means 
of the solution of a Dirichlet problem for the Laplace equation. We gave necessary 
and sufficient conditions on the datum which generates a 4-point which suggest 
that the most probable configurations in nature are those with 3-points.  

\subsection*{Acknowledgment}
The authors thank professor Andrea Dall'Aglio for many
fruitful discussions. The authors are also very grateful 
to the anonymous referees for their insightful comments 
and helpful suggestions.

\end{document}